\documentclass[bjps]{imsart}

\RequirePackage[OT1]{fontenc}
\usepackage{amsthm,amsmath,natbib}
\usepackage{amssymb,bm}
\usepackage{bigints}
\usepackage{float}
\usepackage{graphicx}

\arxiv{arXiv:1803.11271}
\newtheorem{theorem}{Theorem}

\newtheorem{remark}{Remark}
\newtheorem{corollary}{Corollary}
\newtheorem{definition}{Definition}

\newtheorem{asm}{Assumption}
\newenvironment{asmbis}[1]
{\renewcommand{\theasm}{\ref{#1}$'$}%
	\addtocounter{asm}{-1}%
	\begin{asm}}
	{\end{asm}}
\newtheorem{ex}{Example}
\newenvironment{exbis}[1]
{\renewcommand{\theex}{\ref{#1}$'$}%
	\addtocounter{ex}{-1}%
	\begin{ex}}
	{\end{ex}}

\startlocaldefs
\numberwithin{equation}{section}
\theoremstyle{plain}

\endlocaldefs

\begin{document}

\begin{frontmatter}
\title{Reduction principle for functionals of strong-weak dependent vector random fields}
\thankstext{t1}{\textbf{Supplementary Materials} R code used for simulations in this article is available in the folder ``Research materials'' from \href{https://sites.google.com/site/olenkoandriy/}.}

\begin{aug}
\author{\fnms{Andriy} \snm{Olenko}\thanksref{a}\ead[label=e1]{a.olenko@latrobe.edu.au}} and
\author{\fnms{Dareen} \snm{Omari}\thanksref{a}\ead[label=e2]{omari.d@students.latrobe.edu.au}}


\runauthor{A. Olenko \and D. Omari}

\affiliation[a]{La Trobe University}

\address{Department of Mathematics and Statistics\\ La Trobe University\\
Melbourne, VIC, 3086, Australia\\
\printead{e1}\\
\printead{e2}}


\end{aug}

\begin{abstract}
We prove the reduction principle for asymptotics of functionals of vector random fields with weakly and strongly dependent components.
These functionals can be used to construct new classes of random fields with skewed and heavy-tailed distributions. Contrary to the case of scalar long-range dependent random fields, it is shown that the asymptotic behaviour of such functionals is not necessarily determined by the terms at their Hermite rank. 
The results are illustrated by an application to the first Minkowski functional of the Student random fields.
Some simulation studies based on the theoretical findings are also presented.

\end{abstract}




\end{frontmatter}

\section{Introduction}
In various applications researchers often encounter cases involving dependent observations over time or space.
Dependence properties of a random process are usually characterized by the asymptotic behaviour of its covariance function. In particular, a stationary random process $\eta_1(x), \ x\in \mathbb R$, is called weakly (short-range) dependent if its covariance function $B(x)=\textbf{Cov}(\eta_1(x+y),\eta_1(y))$ is integrable, i.e. $\int_{\mathbb{R}}|B(x)|dx<\infty$. On the other hand,  $\eta_1(x)$  possesses strong (long-range) dependence if its covariance function decays slowly and is non-integrable. An alternative definition of long-range dependence is based on singular properties of the spectral density of a random process, such as unboundedness at zero (see~\cite{doukhan2002theory,Souza2008, leonenko2013tauberian}).

Long-range dependent processes play a significant role in a wide range of areas, including finance, geophysics, astronomy, hydrology, climate and engineering (see~\cite{leonenko1999limit,ivanov1989statistical,doukhan2002theory}). In statistical applications, long-range dependent models require developing new statistical methodologies, limit theorems and parameter estimates compared to the weakly dependent case (see \cite{ivanov1989statistical,worsley1994local,leonenko2013tauberian, beran2013long}).

In statistical inference of random fields, limit theorems are the central topic. These theorems play a crucial role in developing asymptotic tests in the large sample theory. The central limit theorem (CLT) holds under the classical normalisation $n^{-d/2}$ when the summands or integrands are weakly dependent random processes or fields. This result was proved by~\cite{breuer1983central} for nonlinear functionals of Gaussian random fields.
A generalisation for stationary Gaussian vector processes was obtained in~\cite{de1995central}, for integral functionals of Gaussian processes or fields in~\cite{chambers1989central},~\cite{hariz2002limit},~\cite{leonenko2014sojourn}, for quasi-associated random fields under various conditions in~\cite{bulinski2012central} and \cite{demichev2015functional}.
Some other CLTs for functionals of Gaussian processes or fields can be found in~\cite{doukhan1999new,coulon2000triangular} and \cite{kratz2017central}. 

The non-central limit theorems arise in the presence of long-range dependence. They use normalising coefficients different than in the CLT and have non-Gaussian limits. These limits are known as Hermite type distributions. 
A non-Gaussian asymptotic was first obtained  in~\cite{rosenblatt1961independence} as a limit for quadratic functionals of stationary Gaussian sequences.   
The article~\cite{taqqu1975weak} continued this research and investigated weak limits of partial sums of Gaussian processes using characteristic functions. The Hermite processes of the first two orders were used.
Later on,~\cite{dobrushin1979non} and~\cite{taqqu1979convergence} established pioneering results in which asymptotics were presented in terms of multiple Wiener-It\^{o} stochastic integrals. 
A generalisation for stationary Gaussian sequences of vectors was obtained in~\cite{arcones1994limit} and~\cite{major2019non}. Multivariate limit theorems for functionals of stationary Gaussian series were addressed under long-range dependence, short-range dependence and a mixture of both in~\cite{bai2013multivariate}. The asymtotices for Minkowski functionals of stationary and isotropic Gaussian random fields with dependent structures were studied in~\cite{ivanov1989statistical}.
\cite{leonenko2014sojourn} obtained the limit theorems for sojourn measures of heavy-tailed random fields (Student and Fisher-Snedecor) under short or long-range dependence assumptions.
Excellent surveys of limit theorems for shortly and strongly dependent random fields can be found in       ~\cite{anh2015rate,doukhan2002theory,ivanov1989statistical,leonenko1999limit,spodarev2014limit}.

The reduction theorems play an important role in studying the asymptotics for random processes and fields.
These theorems show that the asymptotic distributions for functionals of random processes or fields coincide with distributions of other functionals that are much simpler and easier to analyse.  
The CLT can be considered as the ``extreme'' reduction case, when, due to weak dependence and despite the type of functionals and components distributions, asymptotics are reduced to the Gaussian behaviour. The classical non-central limit theorems are based on another ``proper'' reduction principle, when the asymptotic behaviour is reduced only to the leading Hermite term of nonlinear functionals. 
Recently,~\cite{olenko2018reduction} proved the reduction principle for functionals of strongly dependent vector random fields. Components of such vector fields can possess different long-range dependences. 
It was shown that, in contrast to the scalar cases, the limits can be degenerated or can include not all leading Hermite terms.

The available literature, except a few publications, addresses limit theorems and reduction principles for functionals of weakly or strongly dependent random fields separately. 
For scalar-valued random fields it is sufficient as such fields can exhibit only one type of dependence. However, for vector random fields there are various cases with  different dependence structures of components. Such scenarios are important when one aggregates spatial data with different properties. For example, brain images of different patients or GIS data from different regions. Another reason for studying such models is constructing scalar random fields by a nonlinear transformation of a vector field. This approach was used to obtain non-Gaussian fields with some desirable properties, for example, skewed or heavy tailed marginal distributions, see Example~\ref{pp}, Theorem~\ref{th5} and~\cite{leonenko2014sojourn}.

This paper considers functionals of vector random fields which have both strongly and weakly dependent components.
The results in the literature dealt with cases where the interplay between terms at the Hermite rank level and the memory parameter (covariance decay rate) of a Gaussian field completely determines the asymptotic behavior.
This paper shows that in more general settings terms at non-Hermite rank levels can interplay with the memory parameter to determine the limit.
As an application of the new reduction principle we provide some limit theorems for vector random fields. In particular, we show that it is possible to obtain non-Gaussian behaviour for the first Minkowski functional of the Student random field built on different memory type components. It contrasts to the known results about the cases of same type memory components in~\cite{leonenko2014sojourn} where, despite short or long range dependence, only the Gaussian limit is possible.           

The remainder of the paper is organised as follows.
In Section~\ref{sec2} we outline basic notations and definitions that are required in the subsequent sections. Section~\ref{sec3} presents assumptions and main results for functionals of vector random fields with strongly and weakly dependent components.
Sections~\ref{sec44} gives the proofs. 
Section~\ref{sec5} demonstrates some numerical studies.
Short conclusions and some new problems are presented in Section~\ref{sec6}.

\section{Notations}\label{sec2}
This section presents basic notations and definitions of the random field theory and multidimensional Hermite expansions.
Also, we introduce the definition and basic properties of the first Minkowski functional (see~\cite{adler2009random}).
Denote by $\vert\cdot\vert$ and $\Vert\cdot\Vert$ the Lebesgue measure and the Euclidean distance in $\mathbb{R}^{d}$, respectively.	
The symbol $C$ denotes constants that are not important for our exposition. Moreover, the same symbol may be used for different constants appearing in the same proof. 
We assume that all random fields are defined on the same probability space $\left(\Omega,\mathcal{F},\textbf{P}\right)$.
\begin{definition}\label{def3}{\rm\cite{bingham1989regular}} A measurable function ${L}:(0,\infty)\rightarrow (0,\infty)$ is slowly varying at infinity if for all $t>0$,
	$\lim_{r\rightarrow \infty} {L}(tr)/{L}(r) =1.$
\end{definition}
A real-valued random field $\eta_1\left(x\right), x \in \mathbb{R}^{d}$, satisfying $\textbf{E}\eta_{1}^{2}(x)<\infty$ is said to be homogeneous and isotropic if its mean function is a constant and the covariance function $B\left(x,y\right)$ depends only on the Euclidean distance $\left\|{x-y}\right\|$ between $x$ and $y$.	

Let $\eta_1 (x)$, $x\in\mathbb{R}^{d}$, be a measurable mean square continuous zero-mean homogeneous isotropic real-valued random field (see~\cite{ivanov1989statistical,leonenko1999limit}) with the covariance function
\begin{align*}
B\left(r\right) :=\textbf{Cov}(\eta_1(x),\eta_1(y))&=\int_{0}^{\infty}2^{(d-2)/2}\Gamma\bigg(\frac{d}{2}\bigg)J_{(d-2)/2}(rz)\\ &\times(rz)^{(2-d)/2} 
d\Phi\left(z\right),\quad   x,\ y \in \mathbb{R}^d,
\end{align*}
where $r=\rVert x-y \lVert$
and $J_{\nu}(\cdot)$ is the Bessel function of the first kind of order $\nu > -1/2$. The finite measure $\Phi\left(\cdot\right)$ is called the isotropic spectral measure of the random field $\eta_1\left(x\right)$, $x\in\mathbb{R}^{d}$. 

The spectrum of the random field $\eta_1(x)$ is absolutely continuous if there exists a function $\varphi(z)$, $z\in[0,\infty)$, such that
$$\Phi(z)=2\pi^{d/2}\Gamma^{-1}(d/2)\int_{0}^{z}u^{d-1}\varphi(u)du,\quad u^{d-1}\varphi(u)\in L_{1}([0,\infty)).$$
The function $\varphi(\cdot)$ is called the isotropic spectral density of the random field~$\eta_1\left(x\right)$.

A random field $\eta_1\left(x\right)$ with an absolutely continuous spectrum has the following isonormal spectral representation
$$
\eta_1\left(x\right)=\int_{\mathbb{R}^{d}}e^{i\langle\lambda,x \rangle}\sqrt{\varphi(\rVert\lambda\lVert)}W(d\lambda),
$$
where $W(\cdot)$ is the complex Gaussian white noise random measure on $\mathbb{R}^{d}$.

Let $\Delta\subset\mathbb{R}^{d}$ be a Jordan-measurable compact connected set  with $|\Delta| > 0$, and $\Delta$ contains the origin in its interior. Also, assume that $\Delta(r)$, $ r > 0$, is the homothetic image of the
set $\Delta$, with the centre of homothety in the origin and the coefficient $r > 0$, that is,
$|\Delta(r)| = r^d |\Delta|$.
\begin{definition}
	The first Minkowski functional is defined as
	$$
	M_{r}\left\lbrace \eta_1\right\rbrace :=|\left\lbrace x\in\Delta(r): \eta_1(x)>a\right\rbrace |=\int_{\Delta(r)}\chi(\eta_1(x)>a)dx,
	$$
	where $\chi(\cdot)$ is the indicator function and $a$ is a constant.
\end{definition}
The functional $M_{r}\left\lbrace \eta_1\right\rbrace$ has a geometrical meaning, namely, the sojourn measure of the random field $\eta_1(x)$.

In the following we will use integrals of the form	$	\int_{\Delta(r)}\int_{\Delta(r)}Q(\Vert x-y \Vert) dxdy$ with various integrable Borel functions $Q(\cdot)$.
Let two independent random vectors $U$ and $V$ in $\mathbb{R}^{d}$ be uniformly distributed inside the set $\Delta(r)$. Consider a function $Q:\mathbb{R} \rightarrow \mathbb{R}$. Then, we have the following representation
\begin{align}\label{eq1}
\int_{\Delta(r)}\int_{\Delta(r)}Q(\Vert x-y \Vert) dxdy &=|\Delta|^{2}r^{2d}\textbf{E}Q(\Vert U-V\Vert)\notag\\
&=|\Delta|^{2}r^{2d}\int_{0}^{diam\lbrace\Delta(r)\rbrace}Q(\rho)\psi_{\Delta(r)}(\rho)d\rho,
\end{align}
where $\psi_{\Delta(r)}(\rho)$, $\rho\geq 0$, denotes the density function of the distance $\Vert U-V\Vert$ between $U$ and $V$.

Using~(\ref{eq1}) for $r=1$ and $Q(\rho)=\frac{1}{\rho^{\alpha_{0}}}$ one obtains for $\alpha_{0}<d$
\begin{align}\label{rem3}
&\int_{\Delta}\int_{\Delta}\dfrac{dxdy}{\Vert x-y\Vert^{\alpha_{0}}}
=\int_{\Delta}\int_{\Delta}\chi\left(\Vert x-y\Vert\leq diam(\Delta)\right)\Vert x-y\Vert^{-\alpha_{0}}dxdy\notag\\
&\leq C\vert\Delta\vert \int_{0}^{diam(\Delta)}\rho^{d-1-\alpha_{0}}d\rho=C\vert\Delta\vert \dfrac{(diam(\Delta))^{d-\alpha_{0}}}{d-\alpha_{0}}<\infty.
\end{align}	
Let $\left(\eta_{1},\ldots,\eta_{2p}\right)$ be a $2p$-dimensional zero-mean Gaussian vector	and $H_{k}(u)$, $k\geq0$, $u\in\mathbb{R}$, be the Hermite polynomials, see~\cite{taqqu1977law}.	

Consider
$$
e_{v}(\omega):= \prod_{j=1}^{p}H_{k_{j}}(\omega_{j}),
$$
where $\omega=(\omega_{1},\dots,  \omega_{p})' \in \mathbb{R}^{p}$, $v=(k_{1},\dots,k_{p})\in \mathbb{Z}^{p},$ and all $k_{j}\geq 0$ for $j = 1,\dots,p.$\

The polynomials $\{e_{v}(\omega)\}_{v}$ form a complete orthogonal system in the Hilbert space 
\[
L_{2}\left(\mathbb{R}^{p},\phi(\lVert\omega\rVert)d\omega\right)=\left \{G: \int_{\mathbb{R}^{p}}G^{2}(\omega)\phi(\lVert\omega\rVert)d\omega<\infty\right\},
\]
where
$$\phi(\lVert\omega\rVert):= \prod_{j=1}^{p}\phi(\omega_{j}),\quad\quad \phi(\omega_{j}):= \dfrac{e^ {-\omega_{j}^2/2}}{\sqrt{2\pi}}. $$	
An arbitrary function $G(\omega)\in L_{2}\left(\mathbb{R}^{p},\phi(\lVert\omega\rVert)d\omega\right)$ admits an expansion with Hermite coefficients $C_{v}$, given as the following:
\begin{align*}\label{22}
G(\omega)=\sum_{k=0}^{\infty}\sum_{v\in N_{k}}\dfrac{C_{v}e_{v}(\omega)}{v!},\quad C_{v}:=\int_{\mathbb{R}^{p}}G(\omega)e_{v}(\omega)\phi(\lVert\omega\rVert)d\omega,
\end{align*}
where $v!= k_{1}!\dots k_{p}!$ and $$N_{k}:=\left\lbrace(k_{1},\dots,k_{p})\in\mathbb{Z}^{p}: \sum_{j=1}^{p}k_{j}=k,k_{j}\geqslant 0,j=1,\dots,p\right\rbrace.$$
\begin{definition} The smallest integer $\kappa\geqslant 1$ such that $C_{v}=0$ for all $v \in N_{j}$, $j=1,\dots,\kappa-1$, but $C_{v}\neq 0$ for some $v \in N_{\kappa}$ is called the Hermite rank of $G(\cdot)$ and is denoted by $H rank G$.
\end{definition}
In this paper, we consider Student random fields which are an example of heavy-tailed random fields. 
To define such fields, we use a vector random field $\bm{\eta}(x)=[\eta_{1}(x),\dots,\eta_{n+1}(x)]^{'}$, $x\in\mathbb{R}^{d}$, with $\textbf{E}\eta_i(x)=0$ 
where $\eta_{i}(x)$, $i=1,\dots, n+1 $, are independent homogeneous isotropic unit variance Gaussian random fields. 

\begin{definition}
	The Student random field {\rm (}t-random field{\rm )} $T_{n}(x),\ x\in\mathbb{R}^{d}$, is defined by
	\begin{equation*}
	T_{n}(x)=\dfrac{\eta_{1}(x)}{\sqrt{{(1/n)(\eta_{2}^{2}(x)+\cdots+\eta_{n+1}^{2}(x))}}},\quad x\in\mathbb{R}^{d}.
	\end{equation*}
\end{definition}
\section{Reduction principles and limit theorems}\label{sec3}
In this section we present some assumptions and the main results. We prove a version of the reduction principle for vector random fields with weakly and strongly dependent components.

In the following we will use the notation $$\bm{\eta}(x)=[\eta_{1}(x),\dots,\eta_{m}(x),\eta_{m+1}(x),\dots,\eta_{m+n}(x)]',\quad x\in\mathbb{R}^{d},$$ for a vector random field with $m+n$ components.
\begin{asm}\label{ass4}
	Let $\bm{\eta}(x)$ be a vector homogeneous isotropic Gaussian random field with independent components, $\textbf{E}\eta(x)=0$ and a covariance matrix ${B}(x)$ such that ${B}(0)=\mathcal{I}$ and
	\begin{displaymath}
	B_{ij}(\|x\|)=\left\{
	\begin{array}{lr}
	0,\quad {\rm if}\quad i\neq j,\\
	\mathcal{I}_{1}\cdot \|x\|^{-\beta}L_{1}\left(\|x\|\right),\ {\rm if}\  i=j=1,\dots,m,\quad \beta > d, \\
	\mathcal{I}_{2} \cdot	\|x\|^{-\alpha}L_{2}\left(\|x\|\right),\ {\rm if}\  i=j=m+1,\dots,m+n,\  \alpha <d/\kappa,
	\end{array}
	\right.
	\end{displaymath} 		
	where $\mathcal{I}$, $\mathcal{I}_{1}$ and $\mathcal{I}_{2}$ are unit matrices of size $m+n$, $m$ and $n$, respectively, $L_{i}\left(\|\cdot\|\right)$, $i=1,2$, are slowly varying functions at infinity.
\end{asm}
\begin{remark}\label{rem2}
	If Assumption~{\rm\ref{ass4}} holds true the diagonal elements of the covariance matrix $B(x)$ are integrable for the first $m$ elements of $\bm{\eta}(x)$, which corresponds to the case of short-range dependence, and non-integrable for the other elements, which corresponds to the case of long-range dependence. For simplicity, this paper investigates only the case of uncorrelated components. 
\end{remark}
\begin{remark}
For $j=m+1,\dots,m+n$ the random field $\bm{\eta}(x)$ in Assumption~{\rm \ref{ass4}} satisfies 
\begin{align}\label{rtrt}
{\bf E}\left(H_{\kappa}(\eta_{j}(x))H_{\kappa}(\eta_{j}(y))\right)=\kappa!B_{jj}^{\kappa}(\Vert x-y\Vert),\quad x,\ y\in\mathbb{R}^{d},
\end{align}
see {\rm \cite{leonenko1999limit}}. Hence, under Assumption~{\rm \ref{ass4}} the right-hand side of {\rm (\ref{rtrt})} is non-integrable when $\alpha<d/\kappa,$ which guarantees the case of long-range dependence.
\end{remark}
Consider the following random variables:
$$K_{r}:=\int_{\Delta(r)}G\left(\bm{\eta}\left(x\right)\right)dx,\quad \quad K_{r,\kappa}:=\sum_{v\in N_{\kappa}}\dfrac{C_{v}}{v!}\int_{\Delta(r)}e_{v}\left(\bm{\eta}\left(x\right)\right)dx,$$
and $$V_{r}:=\sum_{l \geq \kappa+1}\sum_{v\in N_{l}}\dfrac{C_{v}}{v!}\int_{\Delta(r)}e_{v}\left(\bm{\eta}\left(x\right)\right)dx,$$
where $C_{v}(r)$ are the Hermite coefficients and $\kappa$ is the Hermite rank of the function $G(\cdot)$.
Then $$K_r=K_{r,\kappa}+V_r.$$
\begin{remark}
	The random variable $K_r$ is correctly defined, finite with probability 1 and in the mean square sense, see \S3, Chapter IV in {\rm \cite{gihman2004skorokhod}}.  
\end{remark}

We will use the following notations. Consider the set
$$N_{+}:=\{v=(k_{1,l},\dots,k_{m+n,l}):v\in N_{l},\ C_v\neq 0,\ \ l\geq \kappa\}.$$ Let $$\gamma:=\min_{v \in N_{+}}\big(\beta{\sum_{j=1}^{m}{k}_{j,l}}+\alpha\sum_{j=m+1}^{m+n}{k}_{j,l}\big).$$ Note that $\gamma\geq \alpha \kappa$ and there are cases when $\gamma$ can be reached at multiple $v \in N_+$. Therefore, we define the sets
\begin{center}
	$N_{l}^{*}:=\{v=(k_{1,l},\dots,k_{m+n,l}): v\in N_+\cap N_l,\  \beta{\sum_{j=1}^{m}{k}_{j,l}}+\alpha\sum_{j=m+1}^{m+n}{k}_{j,l}=\gamma\}$ and
\end{center}
$$\mathcal{L}_+:=\{l:\ N_l^{*}\neq \varnothing,\ l\geq \kappa \}.$$
Also, we define the random variable
$$K_{r,l}^{*}:=\sum_{v\in N_{l}^{*}}\dfrac{C_{v}}{v!}\int_{\Delta(r)}e_{v}\left(\bm{\eta}\left(x\right)\right)dx.$$
The random variable $K_{r,l}^{*}\not \equiv 0$ if and only if $l \in \mathcal{L}_+$. 

Theorem~1 in~\cite{olenko2018reduction} gives a reduction principle for vector random fields with strongly dependent components. The following result complements it for the case of random fields with strongly and weakly dependent components.
\begin{theorem}\label{the3}
	Suppose that a the vector random field $\bm{\eta}\left(x\right)$, $x\in\mathbb{R}^{d}$, satisfies Assumption~{\rm\ref{ass4}}, $Hrank G(\cdot)=\kappa\geq 1$ and there is at least one $ v=(k_{1,\kappa},\dots,k_{m+n,\kappa})  \in N_\kappa \cap N_+$ such that $\sum_{j=m+1}^{m+n}k_{j,\kappa}=\kappa$. If for $r\rightarrow \infty$ a limit distribution exists for at least one of the random variables
	$$\dfrac{K_{r}}{\sqrt{\textbf{Var}\ (K_{r})}}\quad and\quad \dfrac{K_{r,\kappa}}{\sqrt{\textbf{Var}\ (K_{r,\kappa})}},$$
	then the limit distribution of the other random variable exists as well, and the limit distributions coincide. Moreover, the limit distributions of 	
	$$\dfrac{K_{r,\kappa}}{\sqrt{\textbf{Var}\ (K_{r,\kappa})}}\quad and\quad \dfrac{K_{r,\kappa}^{*}}{
		\sqrt{\textbf{Var}\ (K_{r,\kappa}^{*})}},$$
	are the same.
\end{theorem}
\begin{remark}
	It will be shown in the proof that the assumptions of Theorem~{\rm{\ref{the3}}} guarantee that $\kappa \in \mathcal{L}_+$.
\end{remark}
\begin{remark}\label{rem2}
		It follows from the asymptotic analysis of the variances in Theorem~{\rm\ref{the3}} that $$\textbf{Var}(K_r)\sim \textbf{Var}(K_{r,\kappa})\sim \textbf{Var}(K_{r,\kappa}^{*}), \quad r \to \infty.$$
\end{remark}
\begin{asm}\label{as2}
	Components $\eta_{j}\left(x\right)$, $j=m+1,\dots,m+n$, of $\bm{\eta}(x)$ have the spectral density $f\left(\|\lambda\|\right)$, $\lambda\in \mathbb{R}^{d}$, such that
	$$
	f\left(\|\lambda\|\right)\sim c_{2}\left(d,\alpha\right)\|\lambda\|^{\alpha-d}L_{2}\left(\dfrac{1}{\|\lambda\|}\right), \quad\quad\quad \Vert \lambda\Vert\rightarrow 0,
	$$
	where $$c_{2}\left(d,\alpha\right)=\dfrac{\Gamma\left((d-\alpha)/{2}\right)}{2^{\alpha}\pi^{d/2}\Gamma\left(\alpha/2\right)}.$$
\end{asm}
Denote the Fourier transform of the indicator function of the set $\Delta$ by
$$
\mathcal{K} \left(x\right):=\int_{\Delta}e^{i\langle u,x \rangle}du,\quad x\in\mathbb{R}^{d}.
$$
Let us define the following random variable
\begin{align}\label{rv}
X_{\kappa}:=\int_{\mathbb{R}^{d\kappa}}^{\prime}\mathcal{K}\left(\lambda_{1}+\dots+\lambda_{\kappa}\right)\dfrac{W(d\lambda_{1})\dots W(d\lambda_{\kappa})}{\|\lambda_{1}\|^{(d-\alpha)/2}\dots \|\lambda_{\kappa}\|^{(d-\alpha)/2}},
\end{align}
where $W(\cdot)$ is the Wiener measure on $(\mathbb{R}^{d},\mathcal{B}^{d})$ and $\int_{\mathbb{R}^{d\kappa}}^{\prime}$ denotes the multiple Wiener-It\^{o} integral.
\begin{theorem}\label{th2}
	Let the vector random field $\bm{\eta}(x)$, $x \in \mathbb{R}^{d}$, and the function $G(\cdot)$ satisfy assumptions of Theorem~{\rm\ref{the3}} and Assumption~{\rm\ref{as2}} holds true. Suppose that 	
	$N_{\kappa}^{*}=\{v \in N_{+}: k_{j,\kappa}=\kappa\  for\  some\  j=m+1,\dots,m+n\}.$
	Then, for $r\to \infty$ the random variables
	$$X_\kappa(r):= c_2^{-\kappa/2}(d,\alpha)r^{(\kappa\alpha)/2-d}L_{2}^{-\kappa/2}(r)K_r$$ converge in distribution to the random variable $\sum_{v \in N_{\kappa}^{*}}\dfrac{C_{v}}{\kappa!}X_v$, where $X_{v}$ are independent copies of $X_{\kappa}$ defined by~$\rm(\ref{rv})$.
\end{theorem}

A popular recent approach to model skew distributed random variables is a convolution $Y=\eta_{1}+\tilde{ \eta}_2$, where $\eta_{1}$ is Gaussian and $\tilde{ \eta}_2$ is continuous positive-valued independent random variables. In this case the probability density of $Y$ has the form $f_{Y}(y)=C \phi(y)G(y)$, where $\phi(\cdot)$ is the pdf of $\eta_{1}$ and $G(\cdot)$ is the cdf of $\tilde{ \eta}_2$, which controls the skewness, see \cite{arellano2005fundamental,azzalini2013skew} and \cite{amiri2019bimodal}.
This approach can be extended to the case of random fields as $Y(x)=\eta_{1}(x)+\tilde{ \eta}_2(x)$, $x \in \mathbb{R}^{d}$, resulting in $Y(x)$ with skewed marginal distributions. In the example below we use $\tilde{ \eta}_2(x)=\eta_2^{2}(x)$ and
show that contrary to the reduction principle for strongly dependent vector random fields in \cite{olenko2018reduction} it is not enough to request $HrankG(\cdot)=\kappa$. The assumption of the existence of $v \in N_\kappa \cap N_+$ satisfying $\sum_{j=m+1}^{m+n}k_{j,\kappa}=\kappa$ in Theorem~\ref{the3} is essential.
\begin{ex}\label{pp}
	Let $m=n=1$, $d=2$ and $G(w_1,w_2)=w_1+w_2^{2}-1$. In this case $G(w_1,w_2)=H_1(w_1)+H_2(w_2)$ and $\kappa=1$, but $k_{2,1}=0 \neq \kappa$. So, the assumption of Theorem~{\rm{\ref{the3}}} does not hold and 
	\begin{align*}
	\dfrac{K_{r}}{\sqrt{\textbf{Var}(K_{r})}}\stackrel{D}{\to} c_2(2,\alpha)X_2, \quad r \to \infty,
	\end{align*}
	which is indeed different from the Gaussian limit that is expected for the case $Hrank G=1$.
\end{ex}
To address situations similar to Example~\ref{pp} and investigate wider classes of vector field we introduce the following modification of Assumption~\ref{ass4}.
\begin{asmbis}{ass4}\label{ass44}
	Let $\bm{\eta}(x)$, $x\in\mathbb{R}^{d}$, be a vector homogeneous isotropic Gaussian random field with independent components, $\textbf{E}\bm{\eta}(x)=0$ and a covariance matrix ${B}(x)$ such that ${B}(0)=\mathcal{I}$ and
	\begin{displaymath}
	B_{ij}(\|x\|)=\left\{
	\begin{array}{lr}
	0,\quad{\rm if}\quad i\neq j,\\
	\mathcal{I}_{1}\cdot \|x\|^{-\beta_i}L_{1}\left(\|x\|\right),\ {\rm if}\  i=j=1,\dots,m, \\
	\mathcal{I}_{2} \cdot	\|x\|^{-\alpha_j}L_{2}\left(\|x\|\right),\ {\rm if}\  i=j=m+1,\dots,m+n, 
	\end{array}
	\right.
	\end{displaymath} 		
	where $\beta_i > d$, $i=1,\dots,m$ and $\alpha_j <d$, $j=m+1,\dots,m+n$.
\end{asmbis}
\begin{remark}
	Under Assumption~{\rm\ref{ass44}} the components $\eta_{m+1}(x),\dots,\eta_{m+n}(x)$ are still strongly dependent, but $H_\kappa(\eta_j(x))$, $j=m+1,\dots,m+n$, do not necessarily preserve strong dependence. If $\kappa \alpha_{j}>d$ the Hermite polynomials of $\eta_j(x)$ become weakly dependent.
\end{remark}
The following modifications of $\gamma$, $N_l^{*}$, $\mathcal{L}_+$ and $K_{r,l}^{*}$ will be used to match Assumption~\ref{ass44}:
$$\tilde{\gamma}:=\min_{v \in N_{+}}\bigg({\sum_{j=1}^{m}\beta_j{k}_{j,l}}+\sum_{j=m+1}^{m+n}\alpha_j{k}_{j,l}\bigg),$$


\begin{center}
$\tilde{N}_{l}^{*}:=\{v=(k_{1,l},\dots,k_{m+n,l}): v\in N_+\cap N_l,$ 
	$\phantom{\tilde{N}_{l}^{*}:= \ }  {\sum_{j=1}^{m}\beta_j{k}_{j,l}}+
	\sum_{j=m+1}^{m+n}\alpha_j{k}_{j,l}=\tilde{\gamma}\},$ 
\end{center}
$$\tilde{\mathcal{L}_+}:=\{l:\ \tilde{N}_l^{*}\neq \varnothing,\ l\geq \kappa \},$$ and
$$\tilde{K}_{r,l}^{*}:=\sum_{v\in \tilde{N}_{l}^{*}}\dfrac{C_{v}}{v!}\int_{\Delta(r)}e_{v}\left(\bm{\eta}\left(x\right)\right)dx.$$
In the following we consider only the cases \mbox{${\sum_{j=1}^{m}\beta_j{k}_{j,l}}+\sum_{j=m+1}^{m+n}\alpha_j{k}_{j,l}\neq d.$} The case when the sum equals $d$ requires additional assumptions, see Section~\ref{sec6}, and will be covered in other publications.

Now, we are ready to formulate a generalization of Theorem~\ref{the3}.
\begin{theorem}\label{th3}
	Suppose that a vector random field $\bm{\eta}(x)$, $x \in \mathbb{R}^{d}$, satisfies {Assumption~\rm\ref{ass44}} and $\tilde{\gamma}<d$. If a limit distribution exists for at least one of the random variables 	$$\dfrac{K_{r}}{\sqrt{\textbf{Var}(K_{r})}}\quad and\quad \dfrac{\sum_{l \in \tilde{\mathcal{L}}_+}\tilde{K}_{r,l}^{*}}{\sqrt{\textbf{Var}\big(\sum_{l \in \tilde{\mathcal{L}}_+} \tilde{K}_{r,l}^{*}\big)}},$$
	then the limit distribution of the other random variable exists as well, and the limit distributions coincide when $r \to \infty$.
\end{theorem}
\begin{asmbis}{as2}\label{as22}
	Components $\eta_{j}\left(x\right)$, $j=m+1,\dots,m+n$, of $\bm{\eta}(x)$ have spectral densities $f_j\left(\|\lambda\|\right)$, $\lambda\in \mathbb{R}^{d}$, such that
	$$
	f_j\left(\|\lambda\|\right)\sim c_{2}\left(d,\alpha_j\right)\|\lambda\|^{\alpha_j-d}L_{2}\left(\dfrac{1}{\|\lambda\|}\right), \quad\quad\quad \Vert \lambda\Vert\rightarrow 0.
	$$
\end{asmbis}
\begin{theorem}\label{th4}
	Let Assumption~{\rm\ref{as22}} and conditions of Theorem~{\rm\ref{th3}} hold true. Suppose that $\tilde{N}_{l}^{*}=\{v\in N_+:k_{j_{l},l}=l \ \mbox{for\ some}\ j_{l}=m+1,\dots,m+n 
	\}$ and there exists  a finite or infinite $\lim_{r\rightarrow \infty}L_2(r)$. Then, for $r \to \infty$ the random variables $$\dfrac{K_r}{r^{d-\tilde{\gamma}/2}\sum_{l \in \tilde{\mathcal{L}}_+}L_{2}^{l/2}(r)}$$
	converge in distribution to the random variable 
	\begin{align}\label{++}
	\sum_{l \in \tilde{\mathcal{L}}_+} a_{l} \cdot \sum_{v\in \tilde{N}_{l}^{*}}\dfrac{C_v}{v!} c_{2}^{l/2}(d,\alpha_{j_{l}})X_v,
	\end{align}
	where $X_v$ are independent copies of random variables
	$$\int_{\mathbb{R}^{dl}}^{\prime}\mathcal{K}\left(\lambda_{1}+\dots+\lambda_{l}\right)\dfrac{W(d\lambda_{1})\dots W(d\lambda_{l})}{\|\lambda_{1}\|^{(d-\alpha_{j_{l}})/2}\dots \|\lambda_{l}\|^{(d-\alpha_{j_{l}})/2}},$$
	and the coefficients $a_l$ are finite and defined by
	$a_l:=\lim_{r\rightarrow \infty}\dfrac{L_2^{l/2}(r)}{\sum_{i \in \tilde{\mathcal{L}}_+}L_2^{i/2}(r)}.$
\end{theorem}

\begin{figure}[H]
	\begin{center}
		\centering
		\includegraphics[width=1\linewidth,trim={0 0 0 2cm},clip]{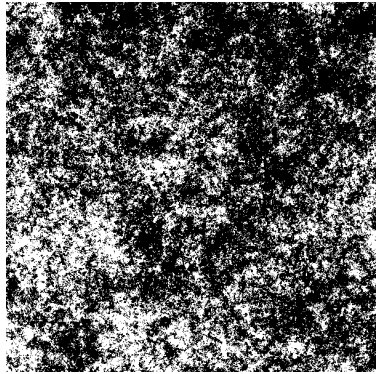}
		\vspace{-2.7cm} 			
		\caption{A two-dimensional excursion set.}  
		\label{fig:ex}	
	\end{center}
\end{figure}

\begin{corollary}
	Let Assumption~{\rm\ref{as22}} and conditions of Theorem~{\rm\ref{th3}} hold true and $n=1$. Then, for $r \to \infty$ the random variable 
	$c_2^{-1}(d,\alpha_{m+1})r^{\tilde{\gamma}/2-d}L_{2}^{-l/2}(r)K_r $	
	converges in distribution to the random variable 
	$(l!)^{-1}C_{(0,\dots,0,l)}~X_{(0,\dots,0,l)},$ where $l=\frac{\tilde{\gamma}}{\alpha_{m+1}}$.
\end{corollary}
\begin{remark}
		It is possible to obtain general versions of Theorems~{\rm\ref{th2}} and~{\rm\ref{th4}} by removing the assumptions about $N_\kappa^{*}$ and $\tilde{ N}_l^{*}$ and requesting only $\mathcal{L}_+\neq \varnothing$ or $\tilde{ \mathcal{L}}_+\neq \varnothing$ respectively. However, it requires an extension of the known non-central limit theorems for vector fields from the discrete to continuous settings, see Section~{\rm\ref{sec6}}. Also, in such general cases the summands in the limit random variables analogous to~{\rm(\ref{++})} would be dependent.
\end{remark}
As an example we consider the first Minkowski functional of Student random fields. The special cases of only weakly or strongly dependent components were studied in~\cite{leonenko2014sojourn}. It was shown that in the  both cases the asymptotic distribution is $N(0,1)$, but with different normalisations, see Theorems~3 and~6 in~\cite{leonenko2014sojourn}.
Figure~\ref{fig:ex} gives a two-dimensional excursion set above the level $a=0.5$ for a realisation of a long-range dependent Cauchy model. 
The excursion set is shown in black colour. 
More details are provided in Section~\ref{sec5}. 
\vspace{-0.2cm}
\vspace{-0.6cm}
The next result shows that for the first Minkowski functional of t-fields obtained from vector random fields with both weakly and strongly dependent components the limit distributions can be non-Gaussian. 
\begin{theorem}\label{th5}
	Let Assumption~{\rm\ref{as22}} hold true, $m=1$, $a\neq 0$, $\alpha_{2}=\dots=\alpha_{m+1}=\alpha<\frac{d}{2}$. Then the random variable
	\[\frac{M_{r}\left\{ T_{n}\right\}-\left| \Delta \right| r^{d}\left(\frac12-\frac12\left(1-I_{\frac{n}{n+a^2}}\left(\frac{n}{2},\frac{1}{2}\right)\right)\cdot {\rm sgn}(a)\right)}{r^{d-\alpha }L(r)}\]
	converges in distribution to the random variable 
	$$\sum_{v\in N_{2}:\ k_{j,2}=2,\ j=2,\dots,n+1}\frac{C_v}{v!}\ c_2(d,\alpha)\ \tilde{X}_v, \quad \mbox{as}\ r \to \infty,$$
		where $\tilde{X}_v$ are independent copies of the random variable
	$$\int_{\mathbb{R}^{2d}}^{\prime}\mathcal{K}\left(\lambda_{1}+\lambda_{2}\right)\dfrac{W(d\lambda_{1}) W(d\lambda_{2})}{\|\lambda_{1}\|^{(d-\alpha)/2} \|\lambda_{2}\|^{(d-\alpha)/2}},$$
	and $sgn(\cdot)$ is the signum function.

\end{theorem}
\begin{remark}
		Random variables $\tilde{ X}_v$ have the Rosenblatt-type distribution, see~{\rm\cite{anh2015rate}}.
\end{remark}
\begin{remark}
		As $m=1$, the first component $\eta_{1}(x)$ is weakly dependent and the remaining components $\eta_{j}(x)$, $j=2,\dots,n+1$, are strongly dependent.
\end{remark}

\section{ Proofs of the results from Section~\ref{sec3}}\label{sec44}
\begin{proof}[\rm\textbf{Proof of Theorem~\ref{the3}}]
	First we study the behaviour of $K_{r,\kappa}$.
	Note, that
	\begin{align*}
	K_{r,\kappa}&=\sum_{v\in N_{\kappa}}\dfrac{C_{v}}{v!}\int_{\Delta(r)}\prod_{j=1}^{m+n}H_{k_{j}}(\eta_{j}(x))dx.
	\end{align*}
	Let us denote the sets $N_{\kappa}^{(i)},\ i=1,2,3,$ as follows
	$$N_{\kappa}^{(1)}:=\{(k_{1,\kappa},\dots,k_{m+n,\kappa}) : \sum_{j=1}^{m}k_{j,\kappa}=\kappa\},$$
	$$N_{\kappa}^{(2)}:=\{(k_{1,\kappa},\dots,k_{m+n,\kappa}): \sum_{j=m+1}^{m+n}k_{j,\kappa}=\kappa\},$$ and
	$$N_{\kappa}^{(3)}:=\{(k_{1,\kappa},\dots,k_{m+n,\kappa}): \sum_{j=1}^{m+n}k_{j,\kappa}=\kappa\  \mbox{and}\  0<\sum_{j=1}^{m}k_{j,\kappa}< \kappa\}.$$
	Then $N_\kappa= \bigcup\limits_{i=1}^{\infty} N_{\kappa}^{(i)}$ and $K_{r,\kappa}$ can be written as  
	\begin{align*}
	K_{r,\kappa}&=\sum_{v_{1}\in N_{\kappa}^{(1)}}\dfrac{C_{v_{1}}}{v_1!}\int_{\Delta(r)}\prod_{j=1}^{m}H_{k_{j,\kappa}}(\eta_{j}(x))dx\\&+
	\sum_{v_2\in N_{\kappa}^{(2)}}\dfrac{C_{v_2}}{v_2!}\int_{\Delta(r)}\prod_{j=m+1}^{m+n}H_{k_{j,\kappa}}(\eta_{j}(x))dx\\
	&+\sum_{v_3\in N_{\kappa}^{(3)}}\dfrac{C_{v_3}}{v_3!}\int_{\Delta(r)}\prod_{j=1}^{m+n}H_{k_{j,\kappa}}(\eta_{j}(x))dx=:\sum_{i=1}^{3}I_i.
	\end{align*}
	Note, that all components $\eta_{j}(x)$, $j=1,\dots,m$, in the first term $I_1$ are weakly dependent and the variance $\textbf{Var}(I_1)$ is equal
	\begin{align*}
	\textbf{Var}(I_1)&=\textbf{Var} \bigg(\sum_{v_1\in N_{\kappa}^{(1)}}\dfrac{C_{v_1}}{v_1!}\int_{\Delta(r)}\prod_{j=1}^{m}H_{k_{j,\kappa}}(\eta_{j}(x))dx\bigg)\\
	&= \sum_{v_1\in N_{\kappa}^{(1)}}\dfrac{C_{v_1}^2}{(v_1!)^2}\int_{\Delta(r)}\int_{\Delta(r)} \textbf{E} \prod_{j=1}^{m}H_{k_{j,\kappa}}(\eta_{j}(x))H_{k_{j,\kappa}}(\eta_{j}(y))\\
	&= \sum_{v_1\in N_{\kappa}^{(1)}}\dfrac{C_{v_1}^2}{v_1!}\int_{\Delta(r)}\int_{\Delta(r)}\prod_{j=1}^{m}B_{jj}^{k_{j,\kappa}}(\Vert x-y \Vert) dxdy\\
	&= \sum_{v_1\in N_{\kappa}^{(1)}}\dfrac{C_{v_1}^2}{v_1!}\int_{\Delta(r)}\int_{\Delta(r)}B_{11}^{\kappa}(\Vert x-y \Vert) dxdy.
	\end{align*}
	Let $u=x-y$ and $v=y$. The Jacobian of this transformation is $\vert J \vert=1$. By denoting $ \Delta(r)-\Delta(r):=\{u \in \mathbb{R}^d:u=x-y$, $x$, $y \in \Delta(r)\}$ then $\textbf{Var}(I_1)$ can be rewritten as
	\[
	\textbf{Var}(I_1)= Cr^d\sum_{v_1\in N_{\kappa}^{(1)}}\dfrac{C_{v_1}^2}{v_1!}\int_{\Delta(r)-\Delta(r)}B_{11}^{\kappa}(\Vert u \Vert) du.
	\]
	It follows from $B_{jj}^{\kappa}(\Vert u \Vert)\leq B_{jj}(\Vert u \Vert)\leq 1$ and by Remark~\ref{rem2} that for weakly dependent components we get 
	\begin{align*}
	\int_{\Delta(r)-\Delta(r)}B_{11}^{\kappa}(\Vert u \Vert)du<\infty \quad \mbox{and} \quad
	\int_{\mathbb{R}^d}B_{11}^{\kappa}(\Vert u \Vert)du<\infty.
	\end{align*}
	Noting that 
	\begin{align*}
	\int_{\Delta(r)-\Delta(r)}B_{jj}^{\kappa}(\Vert u \Vert)du \rightarrow\int_{\mathbb{R}^d}B_{jj}^{\kappa}(\Vert u \Vert)du,\quad r\rightarrow \infty,
	\end{align*}
	one obtains the following asymptotic behaviour of $\textbf{Var}(I_1)$
	\begin{align}\label{wea}
	\textbf{Var}(I_1)\sim Cr^d\sum_{v_1\in N_{\kappa}^{(1)}}\dfrac{C_{v_1}^2}{v_1!}\int_{\mathbb{R}^d}B_{jj}^{\kappa}(\Vert u \Vert) du,\quad r \rightarrow \infty.
	\end{align}
	In contrast, the components $\eta_{j}(x)$, $j=m+1,\dots,m+n$, in the second term $I_2$ are strongly dependent and  
	$\textbf{Var}(I_2)$ can be obtained as follows
	\begin{align*}
	\textbf{Var}(I_2)&=\textbf{Var} \bigg(\sum_{v_2\in N_{\kappa}^{(2)}}\dfrac{C_{v_2}}{v_2!}\int_{\Delta(r)}\prod_{j=m+1}^{m+n}H_{k_{j,\kappa}}(\eta_{j}(x))dx\bigg)\\
	&= \sum_{v_2\in N_{\kappa}^{(2)}}\dfrac{C_{v_2}^{2}}{(v_2!)^{2}}\int_{\Delta(r)}\int_{\Delta(r)} \textbf{E} \prod_{j=m+1}^{m+n}H_{k_{j,\kappa}}(\eta_{j}(x))H_{k_{j,\kappa}}(\eta_{j}(y))\\
	&= \sum_{v_2\in N_{\kappa}^{(2)}}\dfrac{C_{v_2}^{2}}{v_2!}\int_{\Delta(r)}\int_{\Delta(r)}\prod_{j=m+1}^{m+n}[\Vert x-y \Vert ^{-\alpha} L_{2}(\Vert x-y \Vert)]^{k_{j,\kappa}}dxdy\\
	&=r^{2d-\alpha\kappa}\sum_{v_2\in N_{\kappa}^{(2)}}\dfrac{C_{v_2}^{2}}{v_2!}\int_{\Delta}\int_{\Delta}\Vert x-y \Vert ^{-\alpha\kappa} L_{2}^{\kappa}(r\Vert x-y \Vert)dxdy.
	\end{align*}
	By~(\ref{eq1}) we get
	\begin{align*}
	\textbf{Var}(I_2)&=\vert \Delta \vert ^2r^{2d-\alpha\kappa}\sum_{v_2\in N_{\kappa}^{(2)}}\dfrac{C_{v_2}^{2}}{v_2!}\int_{0}^{diam\{\Delta\}}z ^{-\alpha\kappa} L_{2}^{\kappa}(rz)\psi_{\Delta}(z) dz.
	\end{align*}
	Noting that $\alpha \in (0,d/\kappa)$ by Theorem~2.7 in~\cite{seneta1976functions} we obtain
	\begin{align}\label{lon}
	\textbf{Var}(I_2)&\sim c_1(\kappa,\alpha,\Delta)\vert \Delta \vert ^2\sum_{v_2\in N_{\kappa}^{(2)}}\dfrac{C_{v_2}^{2}}{v_2!} r^{2d-\alpha\kappa}L_{2}^{\kappa}(r), \quad r\rightarrow \infty,
	\end{align}
	where $c_1(\kappa,\alpha,\Delta):=\bigintsss_{0}^{diam\{\Delta\}}z ^{-\alpha\kappa}\psi_{\Delta}(z) dz$.
	By~(\ref{rem3}) and the condition $\alpha <d/\kappa$ the coefficient $c_1(\kappa,\alpha,\Delta)$ is finite as 
	\begin{align*}
	c_1(\kappa,\alpha,\Delta)&=\int_{0}^{diam\{\Delta\}} z^{- \alpha\kappa}
	\psi_{\Delta}(z)dz
	=\vert \Delta \vert^{-2} \int_{\Delta}\int_{\Delta} \Vert x-y \Vert^{-\alpha\kappa}dxdy\notag\\
	&\leq \vert \Delta \vert^{-1} \int_{0}^{diam\left\lbrace \Delta\right\rbrace } \rho^{d-\left(1+ \alpha\kappa\right)} d\rho<\infty.
	\end{align*}
	There are strongly and weakly dependent components in the term $I_{3}$ and its variance $\textbf{Var}(I_3)$ can be rewritten as follows
	\begin{align}\label{na}
	\textbf{Var}&(I_3)=\textbf{Var} \bigg(\sum_{v_3\in N_{\kappa}^{(3)}}\dfrac{C_{v_3}}{v_3!}\int_{\Delta(r)}\prod_{j=1}^{m+n}H_{k_{j,\kappa}}(\eta_{j}(x))dx\bigg)\notag\\
	&= \sum_{v_3\in N_{\kappa}^{(3)}}\dfrac{C_{v_3}^{2}}{(v_3!)^{2}}\int_{\Delta(r)}\int_{\Delta(r)} \textbf{E} \prod_{j=1}^{m}H_{k_{j,\kappa}}(\eta_{j}(x))H_{k_{j,\kappa}}(\eta_{j}(y))\notag\\
	&\times\textbf{E}\prod_{j=m+1}^{m+n}H_{k_{j,\kappa}}(\eta_{j}(x))H_{k_{j,\kappa}}(\eta_{j}(y))dxdy\notag\\
	&= \sum_{v_3\in N_{\kappa}^{(3)}}\dfrac{C_{v_3}^{2}}{v_3!}\int_{\Delta(r)}\int_{\Delta(r)}B_{11}^{\sum_{j=1}^{m}k_{j,\kappa}}(\Vert x-y \Vert) B_{m+1m+1}^{\sum_{j=m+1}^{m+n}k_{j,\kappa}}(\Vert x-y \Vert)dxdy\notag\\
	&= \sum_{v_3\in N_{\kappa}^{(3)}}\dfrac{C_{v_3}^{2}}{v_3!}\int_{\Delta(r)}\int_{\Delta(r)} \tilde{B}(\Vert x-y \Vert)dxdy,
	\end{align}
	where 
	\begin{align}\label{cov}
	\tilde{B}(\Vert x-y \Vert)&:= B_{11}^{\sum_{j=1}^{m}k_{j,\kappa}}(\Vert x-y \Vert)  B_{m+1m+1}^{\sum_{j=m+1}^{m+n}k_{j,\kappa}}(\Vert x-y \Vert) \notag \\
	&=\Vert x-y \Vert^{-(\beta{\sum_{j=1}^{m}k_{j,\kappa}}+\alpha\sum_{j=m+1}^{m+n}k_{j,\kappa})}\tilde{L}(\Vert x-y \Vert),
	\end{align}
	and
	\begin{align*}
	\tilde{L}(\Vert x-y \Vert):=L_{1}^{\sum_{j=1}^{m}k_{j,\kappa}}(\Vert x-y \Vert)L_{2}^{\sum_{j=m+1}^{m+n}k_{j,\kappa}}(\Vert x-y \Vert).
	\end{align*}
	Note, that by properties of slowly varying functions $\tilde{L}(\cdot)$ is also a slowly varying function.
	
	If in~(\ref{cov}) the power $\beta{\sum_{j=1}^{m}k_{j,\kappa}}+\alpha\sum_{j=m+1}^{m+n}k_{j,\kappa}$ is greater than $d$ then this case is analogous to the case of $I_{1}$ with short-range dependence and similar to~(\ref{wea}) one obtains
	\begin{align}\label{wee}
	\textbf{Var}(I_3)\sim Cr^d \sum_{v_3\in N_{\kappa}^{(3)}}\dfrac{C_{v_3}^{2}}{v_3!}\int_{\mathbb{R}^d}\tilde{B}(\Vert u \Vert) du,\quad r \rightarrow \infty.
	\end{align}
	This is indeed the case for $N_{\kappa}^{(3)}$ as $\sum_{j=1}^{m}k_{j,\kappa}\geq 1$ and $\beta>d$.

	Note, that by the conditions of the theorem $N_{\kappa}^{(2)}\neq \varnothing$.
	Now, by properties of slowly varying functions (see Preposition~1.3.6 in \cite{bingham1989regular}), we get
	\begin{align}\label{sl}
	\dfrac{\textbf{Var}(I_1)}{\textbf{Var}(I_2)}=\frac{Cr^d\sum_{v_1\in N_{\kappa}^{(1)}}\frac{C_{v_1}^2}{v_1!}
		\bigintsss_{\mathbb{R}^d}B_{11}^{\kappa}(\Vert u \Vert) du}{c_1(\kappa,\alpha,\Delta)\vert \Delta \vert ^2\sum_{v_2\in N_{\kappa}^{(2)}}\frac{C_{v_2}^{2}}{v_2!} r^{2d-\alpha\kappa}L_{2}^{\kappa}(r)}
	\rightarrow 0, \quad r\rightarrow \infty.
	\end{align} 
	By~(\ref{lon}) and~(\ref{wee}) we also obtain
	\begin{align}\label{ml}
	\dfrac{\textbf{Var}(I_3)}{\textbf{Var}(I_2)}&=\dfrac{Cr^d \sum_{v_3\in N_{\kappa}^{(3)}}\frac{C_{v_3}^{2}}{v_3!}\bigintsss_{\mathbb{R}^d}\tilde{B}(\Vert u \Vert) du
	}{c_1(\kappa,\alpha,\Delta)\vert \Delta \vert ^2\sum_{v_2\in N_{\kappa}^{(2)}}\frac{C_{v_2}^{2}}{v_2!} r^{2d-\alpha\kappa}L_{2}^{\kappa}(r)}\rightarrow 0,\ \ r\rightarrow \infty.
	\end{align} 
	
	Note, that 
	\[
	\textbf{Var}\bigg(\sum_{i=1}^{3}I_i\bigg)= \textbf{Var}(I_2)\bigg(\frac{\textbf{Var}(I_1)}{\textbf{Var}(I_2)}+1+\frac{\textbf{Var}(I_3)}{\textbf{Var}(I_2)}+\dfrac{2\sum_{1\leq i <j\leq 3} \textbf{Cov}(I_i,I_j)}{\textbf{Var}(I_2)}\bigg). 
	\]
	 Using the Cauchy–Schwarz inequality 
	$ \vert \textbf{Cov}(I_i,I_j) \vert \leq \sqrt{\textbf{Var}(I_i)\textbf{Var}(I_j)}$ by~(\ref{sl}) and~(\ref{ml}) we get for $r \to \infty$ that 
	\[
	\dfrac{\vert \textbf{Cov}(I_1,I_2) \vert}{\textbf{Var}(I_2)}\leq\sqrt{\dfrac{\textbf{Var}(I_1)}{\textbf{Var}(I_2)}}\rightarrow 0,\quad \dfrac{\vert \textbf{Cov}(I_2,I_3) \vert}{\textbf{Var}(I_2)}\leq\sqrt{\dfrac{\textbf{Var}(I_3)}{\textbf{Var}(I_2)}}\rightarrow 0,\]
\begin{equation}\label{nn}
\dfrac{\vert \textbf{Cov}(I_1,I_3) \vert}{\textbf{Var}(I_2)}\leq\sqrt{\dfrac{\textbf{Var}(I_1)}{\textbf{Var}(I_2)}}\sqrt{\dfrac{\textbf{Var}(I_3)}{\textbf{Var}(I_2)}}\rightarrow 0.
\end{equation}
	Therefore, combining the above results we obtain
	\begin{align}\label{krk}
	\textbf{Var}(K_{r,\kappa})=\textbf{Var}\bigg(\sum_{i=1}^{3}I_i\bigg)\sim \textbf{Var}(I_2)(1+o(1)), \quad r\rightarrow \infty.
	\end{align}
	Now, we study the behaviour of $V_{r}$.
	Similarly to $K_{r,\kappa}$, to investigate   $ \textbf{Var}(V_{r})$ we define the following sets $$N_{l}^{(1)}=\{(k_{1,l},\dots,k_{m,l}): \sum_{j=1}^{m}k_{j,l}=l\},$$ $$N_{l}^{(2)}=\{(k_{1,l},\dots,k_{m+n,l}): \sum_{j=m+1}^{m+n}k_{j,l}=l\},$$ and $$N_{l}^{(3)}=\{(k_{1,l},\dots,k_{m+n,l}): \sum_{j=1}^{m+n}k_{j,l}=l \  \mbox{and}\   0<\sum_{j=1}^{m}k_{j,l}< l\}.$$ Then $V_{r}$ can be written as 
	\begin{align*}
	V_{r}&=\sum_{l \geq \kappa+1}\sum_{v_{1}\in N_{l}^{(1)}}\dfrac{C_{v_{1}}}{v_1!}\int_{\Delta(r)}\prod_{j=1}^{m}H_{k_{j,l}}(\eta_{j}(x))dx\\&+
	\sum_{l \geq \kappa+1}\sum_{v_2\in N_{l}^{(2)}}\dfrac{C_{v_2}}{v_2!}\int_{\Delta(r)}\prod_{j=m+1}^{m+n}H_{k_{j,l}}(\eta_{j}(x))dx\\
	&+\sum_{l \geq \kappa+1}\sum_{v_3\in N_{l}^{(3)}}\dfrac{C_{v_3}}{v_3!}\int_{\Delta(r)}\prod_{j=1}^{m+n}H_{k_{j,l}}(\eta_{j}(x))dx=:\sum_{i=1}^{3}I_i^{(l)}.
	\end{align*}
	Hence, $$
	\textbf{Var}(V_{r}) =\ \textbf{Var}\bigg(\sum_{i=1}^{3}I_i^{(l)}\bigg)=\sum_{i=1}^{3}\textbf{Var}(I_i^{(l)})+2\sum_{1\leq i <j\leq 3} \textbf{Cov}(I_i^{(l)},I_j^{(l)}).
	$$
	The components $\eta_{j}(x)$, $j=1,\dots,m$, in $I_1^{(l)}$ are weakly dependent and $\textbf{Var}(I_1^{(l)})$ is given by
	\raggedbottom
	\begin{align}
	\textbf{Var}(I_1^{(l)})&=\sum_{l \geq \kappa+1}\sum_{v_1\in N_{l}^{(1)}}\dfrac{C_{v_1}^2}{v_1!}\int_{\Delta(r)}\int_{\Delta(r)}\prod_{j=1}^{m} B_{jj}^{k_{j,l}}(\Vert x-y \Vert) dxdy \notag\\
    \label{eq5} &=\sum_{l \geq \kappa+1}\sum_{v_1\in N_{l}^{(1)}}\dfrac{C_{v_1}^2}{v_1!}\int_{\Delta(r)}\int_{\Delta(r)} B_{jj}^{l}(\Vert x-y \Vert) dxdy.
	\end{align}
	As $B_{jj}(\cdot)\leq1$ and $l> \kappa$ we can estimate the expression in~(\ref{eq5}) by
	\begin{align*}
	\textbf{Var} \left( I_1^{(l)}\right) 
	& \leq
	\sum_{l \geq \kappa+1}\sum_{v_1\in N_{l}^{(1)}}\dfrac{C_{v_1}^2}{v_1!}\int_{\Delta(r)}\int_{\Delta(r)} B_{jj}^{\kappa}(\Vert x-y \Vert) dxdy.
	\end{align*}
	It follows from this estimates and the asymptotic~(\ref{wea}) for $\textbf{Var} \left( I_1\right)$ that
	\begin{align}\label{star}
	\textbf{Var} \left( I_1^{(l)}\right)&\leq Cr^{d}\sum_{l \geq \kappa+1}\sum_{v_1\in N_{l}^{(1)}}\dfrac{C_{v_1}^2}{v_1!}\int_{\mathbb{R}^d}B_{11}^{\kappa}(\Vert u \Vert) du, \quad r\rightarrow \infty.
	\end{align}
	In the term $I_{3}^{(l)}$ the components  are strongly and weakly dependent random fields. Similarly to the case of $\textbf{Var}\big(I_3\big)$  we obtain that $\textbf{Var}\big(I_3^{(l)}\big)$ is equal	
	\begin{align*}&\textbf{Var} \bigg(\sum_{l \geq \kappa+1}\sum_{v_3\in N_{l}^{(3)}}\dfrac{C_{v_3}}{v_3!}\int_{\Delta(r)}\prod_{j=1}^{m}H_{k_{j,l}}(\eta_{j}(x))\prod_{j=m+1}^{m+n}H_{k_{j,l}}(\eta_{j}(x))dx\bigg)\\		
		&= \sum_{l \geq \kappa+1}\sum_{v_3\in N_{l}^{(3)}}\dfrac{C_{v_3}^{2}}{(v_3!)^{2}}\int_{\Delta(r)}\int_{\Delta(r)} \textbf{E}\bigg( \prod_{j=1}^{m}H_{k_{j,l}}(\eta_{j}(x))H_{k_{j,l}}(\eta_{j}(y))\bigg)\\
&\times\textbf{E}\bigg(\prod_{j=m+1}^{m+n}H_{k_{j,l}}(\eta_{j}(x))H_{k_{j,l}}(\eta_{j}(y))\bigg)dxdy\\
&=\sum_{l \geq \kappa+1}\sum_{v_3\in N_{l}^{(3)}}\dfrac{C_{v_3}^{2}}{v_3!}\int_{\Delta(r)}\int_{\Delta(r)}B_{11}^{\sum_{j=1}^{m}k_{j,l}}(\Vert x-y \Vert) B_{m+1m+1}^{\sum_{j=m+1}^{m+n}k_{j,l}}(\Vert x-y \Vert)dxdy\\
	&=\sum_{l \geq \kappa+1}\sum_{v_3\in N_{l}^{(3)}}\dfrac{C_{v_3}^{2}}{v_3!}\int_{\Delta(r)}\int_{\Delta(r)}\hat{B}(\Vert x-y \Vert)dxdy,	
	\end{align*}
	where 
	\begin{align*}
	\hat{B}(\Vert x-y \Vert):= \Vert x-y \Vert^{-(\beta{\sum_{j=1}^{m}{k}_{j,l}}+\alpha\sum_{j=m+1}^{m+n}{k}_{j,l})}\hat{L}(\Vert x-y \Vert),
	\end{align*}
	and
	\begin{align*}
	\hat{L}(\Vert x-y \Vert):=L_{1}^{\sum_{j=1}^{m}{k}_{j,l}}(\Vert x-y \Vert)L_{2}^{\sum_{j=m+1}^{m+n}{k}_{j,l}}(\Vert x-y \Vert)
	\end{align*}
	is a slowly varying function.

	Now, as $\sum_{j=1}^{m}k_{j,l}\geq 1$
	then $\beta{\sum_{j=1}^{m}{k}_{j,l}}+\alpha\sum_{j=m+1}^{m+n}{k}_{j,l}>d$ and similar to~(\ref{wea}) the variance $\textbf{Var}\big(I_3^{(l)}\big)$ has the asymptotic behaviour
	\begin{align}\label{lweak}
	\textbf{Var}\big(I_3^{(l)}\big)\sim Cr^d \sum_{l \geq \kappa+1}\sum_{v_3\in N_{l}^{(3)}}\dfrac{C_{v_3}^{2}}{v_3!}\int_{\mathbb{R}^d}\hat{B}(\Vert u \Vert) du,\quad r \rightarrow \infty.
	\end{align}
	Note, that all assumptions of Theorem~1 in~\cite{olenko2018reduction} are satisfied in our case as $\alpha_j=\alpha$, $j=1,\dots,m$, and then 
	$$\sum_{j=1}^{m}\alpha_{j}k_{j,\kappa}=\alpha\kappa\leq (\kappa+1)\min_{1\leq j \leq m}\alpha_{j}=(\kappa+1)\alpha.$$ 
	Therefore, by Theorem~1 in~\cite{olenko2018reduction} we get 
	\begin{align}\label{star2}
	\dfrac{\textbf{Var}(I_2^{(l)})}{\textbf{Var}(I_2)}\to 0,\quad r \to \infty.	
	\end{align}	
	Finally, combining~(\ref{star}),~(\ref{star2})~(\ref{lweak}) and applying the Cauchy–Schwarz inequality analogously to $\textbf{Var}(K_{r,\kappa})$ one obtains that $
	\dfrac{\textbf{Var}(V_r)}{\textbf{Var}(K_{r,\kappa})}\to 0,\ r \to \infty$,
	which proves the asymptotic equivalence of 
	$\dfrac{K_{r}}{\sqrt{Var(K_{r})}}$ and $ \dfrac{K_{r,\kappa}}{\sqrt{Var(K_{r,\kappa})}}.$
	
	It follows from Assumption~\ref{ass4} that
	$$\beta{\sum_{j=1}^{m}{k}_{j,l}}+\alpha\sum_{j=m+1}^{m+n}{k}_{j,l}>\alpha\sum_{j=m+1}^{m+n}k_{j,\kappa}=\alpha\kappa$$
	for all $v=(k_{1,l},\dots,k_{m+n,l})\in N_+ \backslash N_\kappa^{(2)}$ and any $v_2=(0,\dots,0,k_{m+1,\kappa},\dots,$ $ k_{m+n,\kappa})\in N_\kappa^{(2)}$. Hence, $\gamma=\alpha\kappa$, $N_{\kappa}^{*}=N_\kappa^{(2)}\cap N_+=\{(0,\dots,0,k_{m+1,\kappa},\dots,$ $k_{m+n,\kappa})\in N_+\}\neq \varnothing$ and $\mathcal{L}_+=\{\kappa\}$.
	For $v_2\in N_\kappa^{(2)}$ the coefficient $C_{v_{2}}\neq 0$ only if $v_{2}\in N_\kappa^{*}$. Thus, by~(\ref{krk}) we obtain that
	$\dfrac{K_{r,\kappa}}{\sqrt{\textbf{Var}(K_{r,\kappa})}}$ and $\dfrac{K_{r,\kappa}^{*}}{\sqrt{\textbf{Var}(K_{r,\kappa}^{*})}},$
	have the same limit distribution, which completes the proof.
\end{proof}

\begin{proof}[\rm\textbf{Proof of Theorem~\ref{th2}}]
	By Theorem~\ref{the3} 
	$$X_{\kappa}(r)=\dfrac{\sqrt{\textbf{Var}(K_{r})}}{c_2^{\kappa/2}(d,\alpha)r^{d-(\kappa\alpha)/2}L_{2}^{\kappa/2}(r)}\cdot \dfrac{K_r}{\sqrt{\textbf{Var}(K_{r})}}$$
	and
	$$X_{\kappa}^{*}(r):=\dfrac{\sqrt{\textbf{Var}(K_{r})}}{c_2^{\kappa/2}(d,\alpha)r^{d-(\kappa\alpha)/2}L_{2}^{\kappa/2}(r)}\cdot 
	\dfrac{K_{r,\kappa}^{*}}{\sqrt{\textbf{Var}(K_{r,\kappa}^{*})}}$$
	have the same limit distribution if it exists.
	
	By Remark~\ref{rem2}, 
	$\sqrt{\textbf{Var}(K_r)}\sim \sqrt{\textbf{Var}(K_{r,\kappa}^{*})}$, $r \to \infty$, and hence $X_\kappa^{*}(r)$ and $c_2^{-\kappa/2}(d,\alpha)r^{(\kappa\alpha)/2-d}L_{2}^{-\kappa/2}K_{r,\kappa}^{*}$ have the same limit distribution.
	
	$K_{r,k}^{*}$ is a sum of independent terms of the form 
	$$\frac{C_v}{v!}\int_{\Delta(r)}H_{\kappa}(\eta_{j}(x))dx, \quad j=m+1,\dots,m+n.$$

	It follows from Theorem~5 in \cite{leonenko2014sojourn} that for the independent components $\eta_{j}(x)$ and for each $v \in N_\kappa^{*}$
	$$c_2^{-\kappa/2}(d,\alpha)r^{(\kappa\alpha)/2-d}L_{2}^{-\kappa/2}(r)\int_{\Delta(r)}H_{\kappa}(\eta_{j}(x))dx\to X_\kappa,\quad r\to \infty,$$
	which completes the proof.
\end{proof}

\begin{proof}[\rm\textbf{Proof of Example~\ref{pp}}]
	
	It follows from the form of $G(w_1,w_2)$ that  $$K_{r,1}=\int_{\Delta(r)}\eta_{1}(x)dx\quad\ \mbox{and}\quad\  V_{r}=\int_{\Delta(r)}\eta_{2}^{2}(x)dx-\Delta(r).$$ Then by Theorem~1 in~\cite{leonenko2014sojourn}
	\begin{align}\label{pl}
	\dfrac{K_{r,1}}{\sqrt{\textbf{Var}(K_{r,1})}}\stackrel{D}{\to} N(0,1), \quad r \to \infty,
	\end{align}
	and by Theorem~5 in~\cite{leonenko2014sojourn}
	\begin{align*}
	\dfrac{V_{r}}{\sqrt{\textbf{Var}(V_{r})}}\stackrel{D}{\to}c_2(2,\alpha)X_2, \quad r \to \infty.
	\end{align*}
	Using the independence of $\eta_1(\cdot)$ and $\eta_2(\cdot)$,~(\ref{wea}) and applying~(\ref{lon}) to $H_2(\cdot)$ one obtains 
	\begin{align*}
	\textbf{Var}(K_r)=\textbf{Var}(K_{r,1})+\textbf{Var}(V_r)\sim C_1 r^2+C_2r^{2(2-\alpha)}L_2^{2}(r), \quad r \to \infty.
	\end{align*}
	Therefore, for $\alpha \in (0,1)$
	$$\frac{\textbf{Var}(K_{r,1})}{r^{2(2-\alpha)}L_2^{2}(r)}\to 0\quad\  \mbox{and}\quad\ \textbf{Var}(K_r)\sim C_2 r^{2(2-\alpha)}L_2^{2}(r), \quad r \to \infty. $$
	Hence, $\dfrac{K_{r}}{\sqrt{\textbf{Var}(K_{r})}}\stackrel{D}{\to} c_2(2,\alpha)X_2,\  r \to \infty,$ which is different from the limit distribution in~(\ref{pl}).
\end{proof}

\begin{proof}[\rm\textbf{Proof of Theorem~\ref{th3}}]
	It follows from $\tilde{\gamma}<d$ that there is at least one $v \in N_+$ such that $\sum_{j=m+1}^{m+n}\alpha_j{k}_{j,l}<d$. Moreover, as $\tilde{\gamma}$ can be obtained only for $v \in N_+$ with $\sum_{j=1}^{m}{k}_{j,l}=0$ and $\sum_{j=m+1}^{m+n}{k}_{j,l}=l$ then $\tilde{\mathcal{L}_+}$ is a finite set. Hence, it holds $N_+=N_+^{(1)} \cup N_+^{(2)} \cup N_+^{(3)} $, where 
	$$N_{+}^{(1)}=\{(k_{1,l},\dots,k_{m+n,l}): {\sum_{j=1}^{m}\beta_j{k}_{j,l}}+
	\sum_{j=m+1}^{m+n}\alpha_j{k}_{j,l}> d,\ l\geq \kappa\},$$
	$$N_{+}^{(2)}=\{(k_{1,l},\dots,k_{m+n,l}):
	\sum_{j=m+1}^{m+n}\alpha_j{k}_{j,l}=\tilde{\gamma},
	\ l\geq \kappa\},$$ and 
	$$N_{+}^{(3)}=\{(k_{1,l},\dots,k_{m+n,l}):\tilde{\gamma}< {\sum_{j=1}^{m}\beta_j{k}_{j,l}}+
	\sum_{j=m+1}^{m+n}\alpha_j{k}_{j,l}< d,\ l\geq \kappa\},$$ are disjoint sets.
	
	Using the Hermite expantion of $G(\cdot)$ we obtain
	\begin{align*}
	\bm{K_{r}}&=\sum_{v_{1}\in N_{+}^{(1)}}\dfrac{C_{v_{1}}}{v_1!}\int_{\Delta(r)}\prod_{j=1}^{m+n}H_{k_{j,l}}(\eta_{j}(x))dx\\&+
	\sum_{v_2\in N_{+}^{(2)}}\dfrac{C_{v_2}}{v_2!}\int_{\Delta(r)}\prod_{j=m+1}^{m+n}H_{k_{j,l}}(\eta_{j}(x))dx\\
	&+\sum_{v_3\in N_{+}^{(3)}}\dfrac{C_{v_3}}{v_3!}\int_{\Delta(r)}\prod_{j=1}^{m+n}H_{k_{j,l}}(\eta_{j}(x))dx=:\sum_{i=1}^{3}I_i^{\prime}.
	\end{align*}
	Analogously to~(\ref{na}) and~(\ref{cov}) the variance of each summand in $K_r$ has the form 
	$$C\int_{\Delta(r)}\int_{\Delta(r)} \Vert x-y \Vert^{-({\sum_{j=1}^{m}\beta_j{k}_{j,l}}+
		\sum_{j=m+1}^{m+n}\alpha_j{k}_{j,l})}$$
	$$\times L_{1}^{\sum_{j=1}^{m}k_{j,l}}(\Vert x-y \Vert)L_{2}^{\sum_{j=m+1}^{m+n}k_{j,l}}(\Vert x-y \Vert)dxdy.$$ 
	Then, similarly to~(\ref{wea}) and~(\ref{lon}) we obtain that $\textbf{Var}(I_1^{\prime})\sim C r^{d}$ and $\textbf{Var}(I_2^{\prime})\sim C r^{2d-\tilde{\gamma}}\sum_{l \in \tilde{\mathcal{L}}_+}L_{2}^{l}(r)$, $r \to \infty$, and each term in $I_3^{\prime}$ has the variance that is asymptotically equivalent to  
	$$C r^{2d-({\sum_{j=1}^{m}\beta_j{k}_{j,l}}+
		\sum_{j=m+1}^{m+n}\alpha_j{k}_{j,l})}, \quad r \to \infty.$$
	By the definition of $N_{+}^{(i)}$, $i=1,2,3$, we get 
	$$\frac{\textbf{Var}(I_1^{\prime})}{\textbf{Var}(I_2^{\prime})}\to 0, \quad \quad \frac{\textbf{Var}(I_3^{\prime})}{\textbf{Var}(I_2^{\prime})}\to 0, \quad r \to \infty.$$
	Using the Cauchy-Schwarz inequality analogously to~(\ref{nn}) one obtains 
	$$\textbf{Var} (K_r)=\textbf{Var} \bigg(\sum_{i=1}^{3}I_i^{\prime}\bigg)\sim \textbf{Var}(I_2^{\prime})(1+o(1)),\quad r \to \infty.$$
	Finally, noting that $N_+^{(2)}=
	\bigcup\limits_{l \in \tilde{\mathcal{L}}_+}\tilde{N}_{l}^{*}$ completes the proof.
\end{proof}

\begin{proof}[\rm\textbf{Proof of Theorem~\ref{th4}}]
	By Theorem~\ref{th3}
	$$\dfrac{\sqrt{\textbf{Var}(K_{r})}}{r^{d-\tilde{\gamma}/2}\sum_{l \in \tilde{\mathcal{L}}_+}L_{2}^{l}(r)}\cdot \dfrac{K_r}{\sqrt{\textbf{Var}(K_{r})}}\quad \mbox{and} \quad \dfrac{\sqrt{\textbf{Var}(K_{r})}}{r^{d-\tilde{\gamma}/2}\sum_{l \in \tilde{\mathcal{L}}_+}L_{2}^{l}(r)}\cdot \dfrac{\sum_{l \in \tilde{\mathcal{L}}_+}\tilde{K}_{r,l}^{*}}{\sqrt{\textbf{Var}(K_{r})}}$$
	have the same limit distribution if it exists. It follows from the structure of $\tilde{N}_l^{*}$ that $\sum_{l \in \tilde{\mathcal{L}}_+}\tilde{K}_{r,l}^{*}$ is a sum of terms 
	$$\dfrac{C_v}{v!}\int_{\Delta(r)}H_{l}(\eta_{j_{l}}(x))dx, \quad {j_{l}}=m+1,\dots,m+n.$$
	By Theorem~5 in~\cite{leonenko2014sojourn} for $v \in \tilde{N}_l^{*}$
	\begin{align}\label{vv}
	r^{\tilde{\gamma}/2-d}L_{2}^{-l/2}(r)\int_{\Delta(r)}H_l(\eta_{j_{l}}(x))dx \to c_2^{l/2}(d,\alpha_{j_{l}})X_v.
	\end{align}
	Note, that from $\alpha_{j_{l_{1}}}l_1\neq \alpha_{j_{l_{1}}}l_2$, if $l_1 \neq l_2$, follows that $j_{l_{1}}\neq j_{l_{2}}$, if $l_1,\ l_2 \in \tilde{\mathcal{L}_+}$. Therefore, the term in~(\ref{vv}) are independent for different $j_{l}$.
	
	From the existence of $\lim_{r\rightarrow \infty}L_2(r)$ it follows that 
	\begin{align}\label{vv2}
	\bigg(\sum_{i \in \tilde{\mathcal{L}}_+}L_2^{i/2}(r)\bigg)^{-1}=\frac{L_2^{l/2}(r)}{\sum_{i \in \tilde{\mathcal{L}}_+}L_2^{i/2}(r)}\cdot L_2^{-l/2}(r)\sim a_l L_2^{-l/2}(r),
	\end{align}
	for $l\in \tilde{\mathcal{L}_+}$ and $r \to \infty$.
	
	As $l \in \tilde{\mathcal{L}_+}$ then $L_2^{l/2}(r)\leq \sum_{i \in \tilde{\mathcal{L}}_+}L_2^{i/2}(r)$ and all coefficients $a_l$ are finite.
	
	Finally, by combining~(\ref{vv}), (\ref{vv2}), and noting that $$\sqrt{\textbf{Var}(K_r)}\sim \sqrt{\textbf{Var}\bigg(\sum_{l \in \tilde{\mathcal{L}}_+}\tilde{ K}_{r,l}^{*}\bigg)}, \quad r \to \infty,$$ we obtain the statement of the theorem.
\end{proof}

\begin{proof}[\rm\textbf{Proof of Theorem~\ref{th5}}]
	It was shown in~\cite{leonenko2014sojourn} that \[\int_{\Delta(r)} \left(\chi \left( T_{n}(x)>a\right)-\mathbf{E}\left(\chi \left(T_{n}(x)>a\right)\right)\right)dx=
	\int_{\Delta(r)}{G}\left(\bm{\eta}(x)\right) \, dx,
	\]
	where
	\begin{align*}
{\textstyle	G(w)=\chi \bigg( \frac{w_{1}}{\sqrt{\frac{1}{n}\left(w_{2}^{2}+\cdots
			+w_{n+1}^{2}\right)}}>a\bigg)+\frac12\bigg(1-I_{\frac{n}{n+a^2}}\bigg(\frac{n}{2},\frac{1}{2}\bigg)\bigg)
	 \cdot{\rm sgn}(a)-\frac12}.
	\end{align*}

	Formula~(24) in~\cite{leonenko2014sojourn} gives the Hermite coefficients of $G(w)$ for $v \in N_1$:
	\begin{displaymath}
	C_v=\left\{
	\begin{array}{lr}
	\frac{1}{\sqrt{2\pi}\left(1+a^2/n\right)^{n/2}},\quad\mbox{if}\quad v=(1,0,\cdots,0),\\
	0,\quad \quad \quad \quad \quad \quad\ \mbox{if}\ \quad v \in N_1 \backslash \{(1,0,\cdots,0)\}. 
	\end{array}
	\right.
	\end{displaymath} 	
	Thus, $Hrank G=1$.
	
	As $G(w)$ is an even function of $w_i$, $i=2,\dots,n+1$, then $C_v=0$ for all $v \in N_2$ such that $k_{2,i}=k_{2,j}=1$ for some $i \neq j$. For $v \in N_2$ such that $k_{2,j}=2$ for some $j=2,\dots,n+1$, we obtain 
	\begin{align*}
	C_v&=\int_{\mathbb{R}^{n+1}}{G}(w)H_2(w_j)\phi (\left\| w\right\| )dw\\
	&=\int_{\mathbb{R}^{n+1}}\chi \bigg( \dfrac{w_{1}}{\sqrt{\frac{1}{n}\left(w_{2}^{2}+\dots+w_{n+1}^{2}\right)}}>a\bigg)
	(w_j^{2}-1)\phi (\left\| w\right\| )dw\\
	&+ \bigg(\dfrac12\bigg(1-I_{\frac{n}{n+a^2}}\bigg(\frac{n}{2},\frac{1}{2}\bigg)\bigg) \cdot{\rm sgn}(a)-\dfrac12 \bigg)
	\int_{\mathbb{R}^{n+1}}(w_j^{2}-1)\phi (\left\| w\right\| )dw\\
	&=\dfrac{1}{n}\int_{\mathbb{R}^{n+1}}\chi \bigg( \dfrac{w_{1}}{\sqrt{\frac{1}{n}\left(w_{2}^{2}+\dots+w_{n+1}^{2}\right)}}>a\bigg)
	\bigg(\sum_{j=2}^{n+1} w_j^{2}-n\bigg)\phi (\left\| w\right\| )dw\\
	&=\frac{2\,\pi^{n/2}}{n(2\pi)^{(n+1)/2}\Gamma(n/2)}\int_0^\infty(\rho^{2}-n)\rho^{n-1}\,e^{-\frac{\rho^{2}}{2}}\int_{a\rho/\sqrt{n}}^\infty \,e^{-\frac{w_{1}^{2}}{2}}\,dw_1\,d\rho.
	\end{align*}	
	Now we investigate $C_v$ as a function of $a$:
	\begin{align*}
	\dfrac{d}{da}C_v&=-\frac{\sqrt{\big (1+\frac{a^2}{n}\big)/(2\pi)}}{n^{3/2}\ 2^{(n-2)/2}\ \Gamma(n/2)\sqrt{1+\frac{a^2}{n}}}\cdot
	\int_0^\infty(\rho^{n+2}-n\rho^n)\,e^{-\frac{\rho^{2}}{2}(1+\frac{a^2}{n})}\,d\rho \\
	&=-\frac{1}{n^{3/2}\ 2^{n/2}\ \Gamma(n/2)\sqrt{1+\frac{a^2}{n}}} 
	\big(\mathbf{E}|z|^{n+2}-n\mathbf{E}|z|^{n} \big),
	\end{align*}
	where $z \sim N\bigg(0,\ \dfrac{1}{1+\frac{a^2}{n}}\bigg)$. 
	
	Using the formula for the central absolute moments we obtain 
	\[
	\frac{d}{da}C_v=\frac{\Gamma(\frac{n+1}{2})\big(1-\frac{n+1 }{n+{a^2}}\big)}{\ \sqrt{n\pi }\ \Gamma(n/2)\big({1+\frac{a^2}{n}}\big)^{(n+1)/2}}.
	\]
	Thus, $C_v$ is a strictly increasing function on $(-\infty,-1)\cup(1,\infty)$ and it decreases on $(-1,1)$.
	Note, that $\lim_{a\rightarrow \infty}C_v=0$ and by the formula for the central absolute moments
	$$C_v=C\ \int_0^\infty(\rho^{n+1}-n\rho^{n-1})\,e^{-\frac{\rho^{2}}{2}}\,d\rho=0, \quad \mbox{when}\quad a=0.$$
	Therefore, $C_v \neq 0$ for $a \neq 0$ and $v \in N_2$ such that $k_{2,j}=2$ for some $j=2,\dots,n+1$.
	Hence, we obtain that $\tilde{\mathcal{L}_+}=\{2\}$, $\tilde{N}_2^+=\{v \in N_+ \cap N_2: k_{j,2}=2,\ j=2,\dots ,n+1 \}$ and $a_l=1$. The application of Theorem~\ref{th4} completes the proof.
\end{proof}

\section{Simulation studies}\label{sec5}
In the following numerical examples we use the generalised Cauchy family covariance, see~\cite{gneiting2004stochastic} and~\cite{Schlather2019}, to model components of 
 $\bm{\eta} (x)$, $x\in \mathbb{R}^{2}$.

The Cauchy covariance function is 
\[B(\|x\|)=  (1+\|x\|^{2})^{-\frac{z}{2}},\quad z>0.\]
To simulate long-range dependent components we consider $0<z<1$.
In this range of $z$ the covariance function is non-integrable. For the case of weakly dependent components we use $z>2$
which gives integrable covariance functions.

Limit distributions were investigated using the following procedure. Random fields were simulated on the plane $\mathbb R^2,$ i.e. $d=2,$ using the square observation window $\Delta(r) =\{x\in \mathbb{R
}^{2}:\left| x_i\right| <r,\ i=1,2\}.$ The R software package RandomFields (see~\cite{Schlather2019}) was used to simulate $\eta_{i}(x)$, $x\in \mathbb{R}^{2}$, $i=1,2,3,$ from Cauchy models.
\begin{exbis}{pp}
Here we illustrate the results in Example~{\rm\ref{pp}}.	
The Cauchy model was used to simulate
	$\eta_{i}(x)$, $x\in \mathbb{R}^{2}$, $i=1,2$, satisfying Assumption~{\rm \ref{ass44}}  with $\beta=2.5$ and $\alpha=0.2$ respectively.  $1000$ realisations of $H_1(\eta_1(x))$, $H_2(\eta_2(x))$ and $Y(x)=G(\eta_1(x),\eta_2(x))=H_1(\eta_1(x))+ H_2(\eta_2(x)) $ were generated for the large value $r = 80$ to compute distributions of $K_{r,1}$, $V_r$ and $K_r$ respectively. 
		
	Notice that the random field $Y(x)$ has skewed marginal distributions, see Figure~{\rm \ref{fig:his}}. The coefficient of the skewness equals $1.62$, i.e. the marginal distribution of $Y(x)$ has a heavy right-hand tail.
	
		To compare empirical distributions, Q-Q plots of realisations of $K_r$ versus realisations of $V_r$ and $K_{r,1}$ are produced in Figure~{\rm\ref{fig:fig1}}.
		As large $r$ and the number of realisations were selected for simulations these empirical distributions are close to the corresponding asymptotic distributions.
	\begin{figure}[H]
		\begin{center}
			\centering
			\includegraphics[width=.75\linewidth, height=8cm,clip]{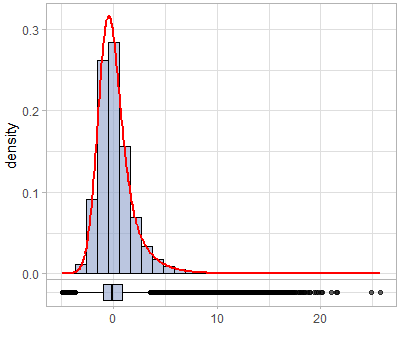}\vspace{-1.cm}
			\caption{The histogram of $Y(x)$.}
			\label{fig:his}
			\vspace{-0.6cm}
		\end{center}		
	\end{figure}
\begin{figure}[H]
	\begin{center}
		\centering
		\includegraphics[width=1\linewidth]{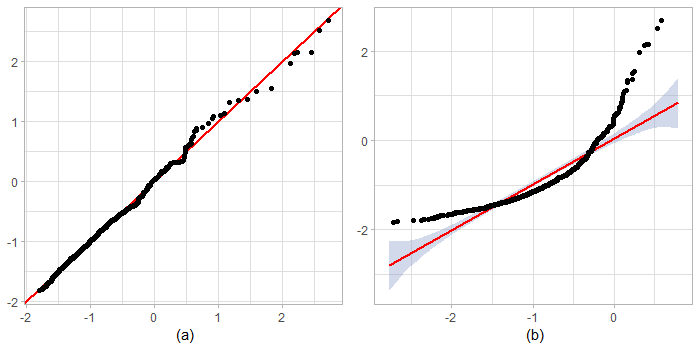} \vspace{-.6cm} 
		\caption{Q-Q plots of {\rm(a)} $K_{r}$ versus $V_r$, {\rm (b)} $K_{r}$ versus $K_{r,1}$}
		\label{fig:fig1}
		\vspace{-0.4cm}
	\end{center}
\end{figure}
 
 It is clear from Figure~{\rm\ref{fig:fig1}(a)} that asymptotic distributions of $K_r$ and $V_r$ are close and the reduction principle works. 
 The Kolmogorov-Smirnov test confirms this result with $p$-value $=0.9937$, see also  Figure~{\rm\ref{fig:ks}} where the plots of empirical cdfs of $K_r$ and $V_r$ are almost identical. 
 However, Figure~{\rm\ref{fig:fig1}(b)} shows that the distributions of $K_r$ and $K_{r,1}$ are different, i.e. asymptotic behaviour of functionals of vector random fields with weak-strong dependent components is not necessarily determined by their Hermite ranks. This result is also confirmed by
 the Kolmogorov-Smirnov $p$-value $= 1.412 \times 10^{-8}$ and Figure~{\rm\ref{fig:ks}}.  	
\begin{figure}[H]
	\begin{center}
		\includegraphics[width=0.8\linewidth,height = 7cm]{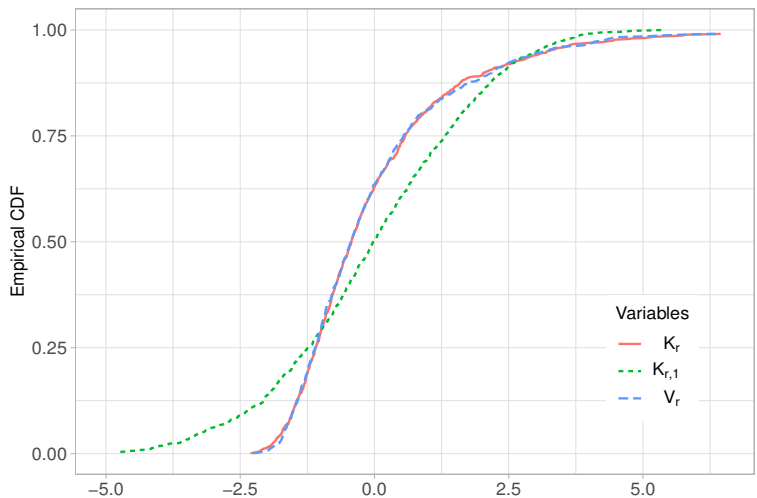} \vspace{-.5cm} 
		\caption{Plots of empirical cdfs of $K_{r}$, $V_r$ and $K_{r,1}$.}
		\label{fig:ks}
		\vspace{-0.4cm}
	\end{center}
\end{figure}
	\end{exbis}
\begin{ex} This example illustrates Theorem~{\rm \ref{th5}}. For $m=1$, $n=2$ and \mbox{$r=600,$} we simulated $500$ realisations of the field $T_{2}(x)=\frac{\eta_{1}(x)}{\sqrt{{(1/2)(\eta_{2}^{2}(x)+\eta_{3}^{2}(x))}}}$, $x \in \mathbb{R}^2.$ For each realisation the area of the excursion set above the level $a = 0.5$ was computed. 
Figure~{\rm \ref{ci}} presents the Q-Q plots with $99\%$ pointwise normal confidence bands for empirical distributions of excursion areas.

The short-range dependent Cauchy model was used to generate realisations of $T_{2}(x)$ with $\beta=4$ for all $\eta_{i}(x)$, $i=1,2,3$. Figure~{\rm \ref{ci}(a)} shows that all the quantiles lie within the confidence bands which  confirms that the first Minkowski functional $M_{r}\{T_{2}(x)\}$ is Gaussian. 

 Another set of realisations of $T_{2}(x)$ was generated using the long-range dependent Cauchy model with $\alpha=0.4$ for all $\eta_{i}(x)$, $i=1,2,3$.  
  Empirical distributions of $M_{r}\{T_{2}(x)\}$ for the obtained realisations of   this long-range dependent and the previous short-range dependent Cauchy models were compared. 
 Figure~{\rm \ref{ci}(b)} shows that the empirical distributions are close and hence, the asymptotic is Gaussian. It is also supported by the  Kolmogorov-Smirnov test with $p$-value $=0.96$.
  Note, that the Gaussianity for these two models follows from the results of Theorems {\rm 3} and {\rm 6} in~{\rm \cite{leonenko2014sojourn}}.
\begin{figure}[H]
	\centering
	\includegraphics[width=1\linewidth]{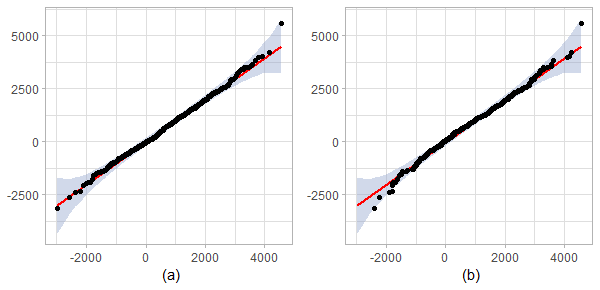} \\
	\includegraphics[width=0.65\linewidth,height=6.5cm]{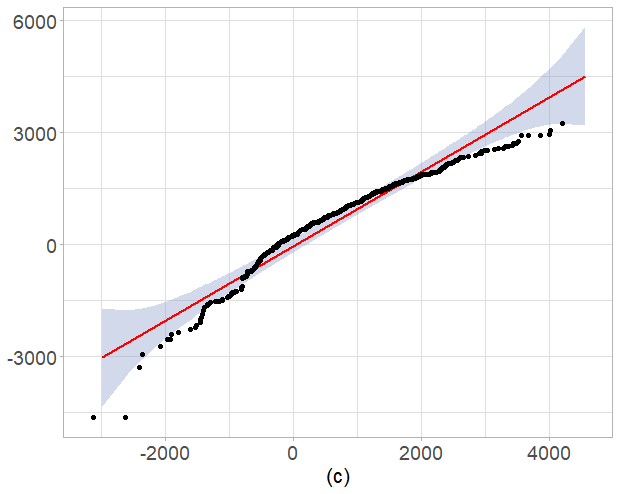}
	\caption{Q--Q plots of realisations of $M_{r}\{T_{2}(x)\}$ of a short-range dependent Cauchy model versus {\rm(a)} the normal distribution, {\rm(b)} $M_{r}\{T_{2}(x)\}$ of a long-range dependent Cauchy model and {\rm(c)} $M_{r}\{T_{2}(x)\}$ of a strong-weak dependent Cauchy model. }
	\label{ci}\vspace{-0.4cm}
\end{figure}

Finally,   $T_{2}(x)$ was generated using the Cauchy fields $\eta_{1}(x)$ with $\beta=4$,  and $\eta_{2}(x)$ and $\eta_{3}(x)$ with $\alpha=0.4$. Note, that $\eta_{1}(x)$ is a weakly dependent component while $\eta_{2}(x)$ and $\eta_{3}(x)$ are strongly dependent ones. The Q-Q plot analogous to Figure~{\rm \ref{ci}(b)} presented in Figure~{\rm \ref{ci}(c)} demonstrates that the distributions are different. 
The corresponding Kolmogorov-Smirnov $p$-value is $0.03$.
In this case the asymptotic distribution of $M_{r}\{T_{2}(x)\}$ of strong-weak dependent components is non-Gaussian and is given by Theorem~{\rm \ref{th5}}.
\end{ex}

\section{Conclusions}\label{sec6}
The paper obtains the reduction principle for vector random fields with strong-weak dependent components. In contrast to the known scalar and vector cases with same type memory components, it is shown that terms at $HrankG $ levels do not necessarily determine limit behaviours. Applications to Minkowski functionals of Student random fields and numerical examples that illustrate the obtained theoretical results are presented. It would be interesting to extend the obtained results to the cases of 
\begin{itemize}
	\item[(1)] cross-correlated components by using some ideas from Theorems~10 and~11 in~\cite{leonenko2014sojourn}. In these theorems it was assumed that the cross-correlation of components is given by some positive definite matrix $\mathcal{A}.$ Then, by using the transformation $\tilde{\eta}=\mathcal{A}^{-1/2}\eta,$ it was possible to transform the vector field to the one with non-correlated components;
	\item[(2)] 	$\sum_{j=1}^{m}\beta_j{k}_{j,l}+\sum_{j=m+1}^{m+n}\alpha_j{k}_{j,l}=d$. It is expected that under some additional assumptions these cases will lead to the CLT, see Remark~2.4 in~\cite{bai2018instability};
	\item[(3)] non-central limit theorems where the condition $k_{j_{l},l}=l \ \mbox{for\ some}\ j_{l}=m+1,\dots,m+n$ is not satisfied. Obtaining such analogous of Theorems~\ref{th4} and \ref{th5}, it requires an extension of Arcones-Major results, see~\cite{arcones1994limit,major2019non}, to continuous settings. While the direct proof may need substantial efforts, see~\cite{major2019non}, one can try the simpler strategy proposed in~\cite{alodat2019asymptotic}. Namely, to prove that discrete and continuous functionals have same limits and then to apply the known discrete result from \cite{arcones1994limit} and \cite{major2019non};
	\item[(4)] cyclically dependent components, i.e. when the spectral density has singular points outside the origin, see, for example, \cite{klykavka2012asymptotic} and \cite{olenko2013limit}. 
\end{itemize}

\section*{Acknowledgements}
Andriy Olenko was partially supported under the Australian Research Council's Discovery Projects funding scheme (project DP160101366). The authors also would like to thank the anonymous referees for their suggestions that helped to improve the paper.

\end{document}